\documentclass[onecolumn]{article}
\usepackage{amsmath,graphicx, amsfonts, amsthm, xfrac, xspace}
\usepackage{lipsum}
\usepackage{algorithm, algpseudocode}
\usepackage{multicol}
\usepackage{amssymb}
\usepackage[caption=false,font=normalsize,labelfont=sf,textfont=sf]{subfig}
\usepackage{url}
\usepackage{graphicx}
\usepackage{cite}
\newcommand{\nuc}{\newcommand}

\nuc{\ISd}{Itakura-Saito divergence}
\nuc{\pp}{point process}
\nuc{\bb}{Bhattacharyya bound}
\nuc{\cvd}{covariance density}
\nuc{\psd}{power spectral density}
\nuc{\BD}{Bregman divergence}
\nuc{\ifunc}{intensity function}

\nuc{\Psp}{Poisson process}
\nuc{\rnp}{renewal process}
\nuc{\Rnp}{Renewal process}
\nuc{\Hp}{Hawkes process}
\nuc{\HL}{Hawkes-Laguerre}
\nuc{\Lt}{Laplace transform}
\nuc{\iLt}{inverse Laplace transform}
\nuc{\iets}{interevent times}
\nuc{\Gd}{Gamma distributed}
\nuc{\lhf}{likelihood function}
\nuc{\lh}{likelihood}
\nuc{\tiv}{time-invariant}
\nuc{\tv}{time-varying}
\nuc{\Bhatt}{Bhattacharyya}
\nuc{\Bnl}{Bernoulli}
\nuc{\cpi}{counting process increment}

\nuc{\cme}{concentrated matrix exponential}
\nuc{\se}{self-exciting}
\nuc{\AWf}{Abate-Whitt framework}

\nuc{\ul}{\underline}

\nuc{\Ra}{\Rightarrow}
\nuc{\La}{\Leftarrow}
\nuc{\LRa}{\Leftrightarrow}
\nuc{\toi}{\to\infty}
\nuc{\trid}{\triangledown}
\nuc{\triu}{\triangleup}

\nuc{\twn}{t_1^n}
\nuc{\xwn}{x_1^n}
\nuc{\swn}{\ssum1^n}
\nuc{\swN}{\ssum1^N}
\nuc{\swP}{\ssum1^P}
\nuc{\swp}{\ssum1^p}
\nuc{\soN}{\ssum0^N}
\nuc{\swm}{\ssum1^m}
\nuc{\swM}{\ssum1^M}
\nuc{\siwn}{\ssum{i=1}^n}
\nuc{\siwnl}{\ssum{i=1}^{n_l}}
\nuc{\srwn}{\ssum{r=1}^n}
\nuc{\srwnl}{\ssum{r=1}^{n_l}}
\nuc{\srwnm}{\ssum{r=1}^{n_m}}
\nuc{\sawM}{\ssum{a=1}^M}
\nuc{\sbwM}{\ssum{b=1}^M}
\nuc{\siwm}{\ssum{i=1}^m}
\nuc{\skwm}{\ssum{k=1}^m}
\nuc{\siwM}{\ssum{i=1}^M}
\nuc{\skwM}{\ssum{k=1}^M}
\nuc{\sjwm}{\ssum{j=1}^m}
\nuc{\sjwM}{\ssum{j=1}^M}
\nuc{\slwM}{\ssum{l=1}^M}
\nuc{\smwM}{\ssum{m=1}^M}
\nuc{\sjwp}{\ssum{j=1}^p}
\nuc{\sjwP}{\ssum{j=1}^P}
\nuc{\siwK}{\ssum{i=1}^K}
\nuc{\swk}{\ssum1^k}
\nuc{\soi}{\ssum0^\infty}
\nuc{\swi}{\ssum1^\infty}
\nuc{\pwn}{\Pi_1^n}
\nuc{\pwk}{\Pi_1^k}
\nuc{\pwm}{\Pi_1^m}
\nuc{\prwn}{\Pi_{r=1}^n}

\nuc{\Stiltr}{\sum_{\ti<\tr}}
\nuc{\Sticlltrcm}{\sum_{\ticl<\trcm}}

\nuc{\Swn}{\sum_1^n}
\nuc{\Swm}{\sum_1^m}
\nuc{\SwM}{\sum_1^M}
\nuc{\Swk}{\sum_1^k}
\nuc{\Pwn}{\prod_1^n}
\nuc{\Pwm}{\prod_1^m}
\nuc{\Swp}{\sum_1^p}
\nuc{\SwP}{\sum_1^P}

\nuc{\intpi}{\int_{-\pi}^{\pi}}
\nuc{\intinfty}{\int_{-\infty}^{\infty}}
\nuc{\intoT}{\int_0^T}
\nuc{\intox}{\int_0^x}
\nuc{\intxi}{\int_x^\infty}
\nuc{\intot}{\int_0^t}
\nuc{\intoinfty}{\int_0^\infty}
\nuc{\intoi}{\int_0^\infty}

\nuc{\limto}{\lim_{t\to0}}
\nuc{\limTo}{\lim_{T\to0}}
\nuc{\limxo}{\lim_{x\to0}}
\nuc{\limso}{\lim_{s\to0}}
\nuc{\stoo}{s\to0}
\nuc{\stoi}{s\to\infty}
\nuc{\xtoi}{x\to\infty}
\nuc{\xtoo}{x\to0}
\nuc{\Ttoi}{T\to\infty}

\nuc{\limti}{\lim_{t\to\infty}}
\nuc{\limTi}{\lim_{T\to\infty}}
\nuc{\limxi}{\lim_{x\to\infty}}
\nuc{\limsi}{\lim_{s\to\infty}}

\nuc{\hlamt}{\hat{\lambda}(t)}

\nuc{\nid}{n_i^\delta}
\nuc{\Nid}{N_i^\delta}

\nuc{\intsum}{\ssum0^\infty\int_{R_n(T)}}
\nuc{\ints}{\ssum0^\infty\int}

\nuc{\wos}{\frac1s}
\nuc{\wox}{\frac1x}
\nuc{\wot}{\frac1t}
\nuc{\woT}{\frac1T}
\nuc{\wo}[1]{\frac1{#1}}

\nuc{\cH}{\mathcal{H}}
\nuc{\cL}{\mathcal{L}}
\nuc{\cK}{\mathcal{K}}
\nuc{\cl}{\mathcal{l}}
\nuc{\cI}{\mathcal{I}}
\nuc{\bbR}{\mathbb{R}}
\nuc{\bbN}{\mathbb{N}}
\nuc{\bx}{{\bf x}}
\nuc{\by}{{\bf y}}

\nuc{\eq}[1]{\begin{align*}#1\end{align*}}
\nuc{\eqn}[1]{\begin{align}#1\end{align}}
\nuc{\bmat}[1]{\begin{bmatrix}#1\end{bmatrix}}
\nuc{\mat}[1]{\begin{matrix}#1\end{matrix}}
\nuc{\smat}[1]{\begin{smallmatrix}#1\end{smallmatrix}}

\nuc{\theo}[1]{\begin{theorem}#1\end{theorem}}
\nuc{\lem}[1]{\begin{lemma}#1\end{lemma}}
\nuc{\defi}[1]{\begin{definition}#1\end{definition}}
\nuc{\exa}[1]{\begin{example}#1\end{example}}
\nuc{\cor}[1]{\begin{corollary}#1\end{corollary}}
\nuc{\prop}[1]{\begin{proposition}#1\end{proposition}}
\nuc{\res}[1]{\begin{result}#1\end{result}}
\nuc{\pro}[1]{\begin{proof}#1\end{proof}}
\nuc{\cas}[1]{\begin{cases}#1\end{cases}}
\nuc{\arr}[2]{\begin{array}{#1}#2\end{array}}
\nuc{\bra}[1]{\left(#1\right)}
\nuc{\sqbra}[1]{\left[#1\right]}
\nuc{\ang}[1]{\langle#1\rangle}
\nuc{\Ver}[1]{\lVert#1\rVert}
\nuc{\ver}[1]{\lvert#1\rvert}
\nuc{\ssum}[1]{{\textstyle\sum}_{#1}}
\nuc{\sprod}[1]{{\textstyle\prod}_{#1}}

\nuc{\sbij}[1]{#1_{ij}}
\nuc{\sbbij}[1]{\bar #1_{ij}}
\nuc{\sbi}[1]{#1_{i}}
\nuc{\sbbi}[1]{\bar #1_{i}}
\nuc{\sbj}[1]{#1_{j}}
\nuc{\sbbj}[1]{\bar #1_{j}}
\nuc{\nfd}[1]{#1^{(n)}}
\nuc{\limo}[1]{\lim_{#1\to0}}
\nuc{\limi}[1]{\lim_{#1\to\infty}}

\nuc{\grad}{\bigtriangledown}

\newcommand{\nlist}[1]{\begin{enumerate}#1\end{enumerate}}

\nuc{\itum}[1]{\item[(#1)]}
\nuc{\ita}{\itum{a}}
\nuc{\itb}{\itum{b}}
\nuc{\itc}{\itum{c}}
\nuc{\itd}{\itum{d}}
\nuc{\ite}{\itum{e}}
\nuc{\itf}{\itum{f}}
\nuc{\itg}{\itum{g}}
\nuc{\ith}{\itum{h}}

\DeclareMathOperator{\E}{E}
\DeclareMathOperator{\var}{var}

\nuc{\ai}{a_i}
\nuc{\ak}{a_k}
\nuc{\aj}{a_j}
\nuc{\aij}{a_{ij}}
\nuc{\akl}{a_{kl}}
\nuc{\akj}{a_{kj}}
\nuc{\abi}{\bar a_i}
\nuc{\abk}{\bar a_k}
\nuc{\abj}{\bar a_j}
\nuc{\abij}{\bar a_{ij}}
\nuc{\abkj}{\bar a_{kj}}
\nuc{\alpi}{\alpha_{i}}
\nuc{\alpk}{\alpha_{k}}
\nuc{\alpj}{\alpha_{j}}
\nuc{\alpij}{\alpha_{ij}}
\nuc{\alpkj}{\alpha_{kj}}
\nuc{\alpab}{\alpha_{ab}}
\nuc{\Ahw}{\hat A_1}
\nuc{\Aho}{\hat A_0}
\nuc{\ajcl}{a_{j,l}}
\nuc{\ajclcm}{a_{j,l,m}}

\nuc{\Bi}{B_i}
\nuc{\Bj}{B_j}
\nuc{\Bij}{B_{ij}}
\nuc{\Bkj}{B_{kj}}
\nuc{\bij}{b_{ij}}
\nuc{\bkj}{b_{kj}}
\nuc{\Bbi}{\bar B_i}
\nuc{\Bbj}{\bar B_j}
\nuc{\Bbij}{\bar B_{ij}}
\nuc{\Bbkj}{\bar B_{kj}}
\nuc{\beti}{\beta_i}
\nuc{\betj}{\beta_j}
\nuc{\betk}{\beta_k}
\nuc{\betaa}{\beta_a}
\nuc{\betb}{\beta_b}
\nuc{\betij}{\beta_{ij}}
\nuc{\betkl}{\beta_{kl}}
\nuc{\betkj}{\beta_{kj}}
\nuc{\betab}{\beta_{ab}}
\nuc{\betjcl}{\beta_{j,l}}
\nuc{\betjclcm}{\beta_{j,l,m}}
\nuc{\Bhw}{\hat B_1}
\nuc{\Bho}{\hat B_0}

\nuc{\ci}{c_i}
\nuc{\cj}{c_j}
\nuc{\ck}{c_k}
\nuc{\ca}{c_a}
\nuc{\cb}{c_b}

\nuc{\cm}{c_m}
\nuc{\cij}{c_{ij}}
\nuc{\ckj}{c_{kj}}
\nuc{\cab}{c_{ab}}
\nuc{\cbi}{\bar c_i}
\nuc{\cbj}{\bar c_j}
\nuc{\Ci}{C_i}
\nuc{\Cj}{C_j}
\nuc{\Cij}{C_{ij}}
\nuc{\Cbi}{\bar C_i}
\nuc{\Cbj}{\bar C_j}
\nuc{\Chw}{\hat C_1}
\nuc{\Cho}{\hat C_0}

\nuc{\Di}{D_i}
\nuc{\Dj}{D_j}
\nuc{\Dbi}{\bar D_i}
\nuc{\Dbj}{\bar D_j}
\nuc{\Dij}{D_{ij}}
\nuc{\Dkj}{D_{kj}}
\nuc{\Dab}{D_{ab}}
\nuc{\dij}{d_{ij}}
\nuc{\dab}{d_{ab}}
\nuc{\dkj}{d_{kj}}
\nuc{\Dbij}{\bar D_{ij}}
\nuc{\Dbkj}{\bar D_{kj}}
\nuc{\Dbab}{\bar D_{ab}}

\nuc{\et}{e_t}

\nuc{\fk}{f_k}
\nuc{\Fi}{F_i}
\nuc{\Fk}{F_k}
\nuc{\Fj}{F_j}
\nuc{\Fa}{F_a}
\nuc{\Fb}{F_b}
\nuc{\Fm}{F_m}
\nuc{\Ft}{F_t}
\nuc{\Fbi}{\bar F_i}
\nuc{\Fbj}{\bar F_j}
\nuc{\Fba}{\bar F_a}
\nuc{\Fbb}{\bar F_b}
\nuc{\Fij}{F_{ij}}
\nuc{\Fkj}{F_{kj}}
\nuc{\Fab}{F_{ab}}
\nuc{\fij}{f_{ij}}
\nuc{\fkj}{f_{kj}}
\nuc{\fab}{f_{ab}}
\nuc{\Fbij}{\bar F_{ij}}
\nuc{\Fbkj}{\bar F_{kj}}
\nuc{\Fbab}{\bar F_{ab}}
\nuc{\Fwtu}{F_{12}}
\nuc{\Fbwtu}{\bar F_{12}}
\nuc{\Fpct}{F_{p,t}}

\nuc{\gam}{\gamma}
\nuc{\Gam}{\Gamma}
\nuc{\gami}{\gam_i}
\nuc{\gamk}{\gam_k}
\nuc{\gamj}{\gam_j}
\nuc{\gamij}{\gam_{ij}}
\nuc{\gamkj}{\gam_{kj}}
\nuc{\gamab}{\gam_{ab}}
\nuc{\Gi}{G_i}
\nuc{\Gk}{G_k}
\nuc{\Gj}{G_j}
\nuc{\Ga}{G_a}
\nuc{\Gb}{G_b}
\nuc{\Gbi}{\bar G_i}
\nuc{\Gbj}{\bar G_j}
\nuc{\Gij}{G_{ij}}
\nuc{\Gkj}{G_{kj}}
\nuc{\Gab}{G_{ab}}
\nuc{\Gbab}{\bar G_{ab}}
\nuc{\Gbij}{\bar G_{ij}}
\nuc{\Gbkj}{\bar G_{kj}}
\nuc{\Gwtu}{G_{12}}
\nuc{\Gbwtu}{\bar G_{12}}

\nuc{\hi}{h_i}
\nuc{\hk}{h_k}
\nuc{\hj}{h_j}
\nuc{\ha}{h_a}
\nuc{\hb}{h_b}
\nuc{\hd}{h_d}
\nuc{\hbi}{\bar h_i}
\nuc{\hbk}{\bar h_k}
\nuc{\hbj}{\bar h_j}
\nuc{\hba}{\bar h_a}
\nuc{\hbb}{\bar h_b}
\nuc{\hij}{h_{ij}}
\nuc{\hbij}{\bar h_{ij}}
\nuc{\hkj}{h_{kj}}
\nuc{\hbkj}{\bar h_{kj}}
\nuc{\hab}{h_{ab}}
\nuc{\hbab}{\bar h_{ab}}
\nuc{\hjcm}{h_{j,m}}
\nuc{\hjcl}{h_{j,l}}
\nuc{\hjclcm}{h_{j,l,m}}
\nuc{\Hi}{H_i}
\nuc{\Hk}{H_k}
\nuc{\Hj}{H_j}
\nuc{\Hbi}{\bar H_i}
\nuc{\Hbk}{\bar H_k}
\nuc{\Hbj}{\bar H_j}
\nuc{\Hij}{H_{ij}}
\nuc{\Hbij}{\bar H_{ij}}
\nuc{\Hkj}{H_{kj}}
\nuc{\Hbkj}{\bar H_{kj}}
\nuc{\Hab}{H_{kj}}
\nuc{\Hbab}{\bar H_{ab}}
\nuc{\Hwtu}{H_{12}}
\nuc{\Hbwtu}{\bar H_{12}}
\nuc{\Hinf}{\cH_\infty}
\nuc{\cHTicj}{\cH_{T,i,j}}

\nuc{\kapij}{\kappa_{ij}}
\nuc{\kapbij}{\bar\kappa_{ij}}
\nuc{\kapkj}{\kappa_{kj}}
\nuc{\kapbkj}{\bar\kappa_{kj}}
\nuc{\kapab}{\kappa_{ab}}
\nuc{\kapbab}{\bar\kappa_{ab}}
\nuc{\kapji}{\kappa_{ji}}
\nuc{\kapbji}{\bar\kappa_{ji}}
\nuc{\Ki}{K_i}
\nuc{\Kk}{K_k}
\nuc{\Kj}{K_j}
\nuc{\Kt}{K_t}
\nuc{\Kbi}{\bar K_i}
\nuc{\Kbk}{\bar K_k}
\nuc{\Kbj}{\bar K_j}
\nuc{\Kij}{K_{ij}}
\nuc{\Kbij}{\bar K_{ij}}
\nuc{\Kkj}{K_{kj}}
\nuc{\Kbkj}{\bar K_{kj}}
\nuc{\Kab}{K_{ab}}
\nuc{\Kbab}{\bar K_{ab}}
\nuc{\Kji}{K_{ji}}
\nuc{\Kbji}{\bar K_{ji}}
\nuc{\Kwtu}{K_{12}}
\nuc{\Kbwtu}{\bar K_{12}}

\nuc{\lam}{\lambda}
\nuc{\Lam}{\Lambda}
\nuc{\lami}{\lambda_{i}}
\nuc{\lamk}{\lambda_{k}}
\nuc{\lamj}{\lambda_{j}}
\nuc{\lama}{\lambda_{a}}
\nuc{\lamb}{\lambda_{b}}
\nuc{\lamm}{\lambda_{m}}
\nuc{\laml}{\lambda_{l}}
\nuc{\lamij}{\lambda_{ij}}
\nuc{\lamkj}{\lambda_{kj}}
\nuc{\lamab}{\lambda_{ab}}
\nuc{\lamicj}{\lambda_{i,j}}
\nuc{\lambij}{\bar\lambda_{ij}}
\nuc{\lambab}{\bar\lambda_{ab}}
\nuc{\Lami}{\Lambda_{i}}
\nuc{\Lamk}{\Lambda_{k}}
\nuc{\Lama}{\Lambda_{a}}
\nuc{\Lamj}{\Lambda_{j}}
\nuc{\Lamij}{\Lambda_{ij}}
\nuc{\Lamkj}{\Lambda_{kj}}
\nuc{\Lamab}{\Lambda_{ab}}
\nuc{\Lamicj}{\Lambda_{i,j}}
\nuc{\Lambij}{\bar\Lambda_{ij}}

\nuc{\lamTt}{\lam_t^T}
\nuc{\lamtgT}{\lam_{t|T}}
\nuc{\lamTu}{\lam_u^T}
\nuc{\lamugT}{\lam_{u|T}}
\nuc{\lamicjclcm}{\lam_{i,j,l,m}}
\nuc{\lamTtcu}{\lam_{t,u}^T}
\nuc{\lamtcugT}{\lam_{t,u|T}}
\nuc{\lamTtpw}{\lam^T_{t+1}}
\nuc{\lamtpwgT}{\lam_{t+1|T}}
\nuc{\LamTt}{\Lam_t^T}
\nuc{\LamtgT}{\Lam_{t|T}}
\nuc{\LamTu}{\Lam_u^T}
\nuc{\LamugT}{\Lam_{u|T}}
\nuc{\LamTtpw}{\Lam^T_{t+1}}
\nuc{\LamtpwgT}{\Lam_{t+1|T}}
\nuc{\LamTtcu}{\Lam_{t,u}^T}
\nuc{\LamtcugT}{\Lam_{t,u|T}}

\nuc{\Li}{L_i}
\nuc{\Lk}{L_k}
\nuc{\Lj}{L_j}
\nuc{\cLa}{\cL_a}
\nuc{\cLb}{\cL_b}

\nuc{\Lb}{L_b}
\nuc{\Lbi}{\bar L_i}
\nuc{\Lbk}{\bar L_k}
\nuc{\Lbj}{\bar L_j}
\nuc{\Lij}{L_{ij}}
\nuc{\Lkj}{L_{kj}}
\nuc{\Lab}{L_{ab}}
\nuc{\Licj}{L_{i,j}}
\nuc{\Lbij}{\bar L_{ij}}
\nuc{\Lbkj}{\bar L_{kj}}
\nuc{\Lwtu}{L_{12}}
\nuc{\Lbwtu}{\bar L_{12}}

\nuc{\mi}{m_i}
\nuc{\mk}{m_k}
\nuc{\mj}{m_j}
\nuc{\mij}{m_{ij}}
\nuc{\mkj}{m_{kj}}
\nuc{\mab}{m_{ab}}
\nuc{\Mij}{M_{ij}}
\nuc{\Mkj}{M_{kj}}
\nuc{\Mab}{M_{ab}}
\nuc{\Mfij}{M_{\fij}}
\nuc{\mui}{\mu_i}
\nuc{\muk}{\mu_k}
\nuc{\muj}{\mu_j}
\nuc{\mua}{\mu_a}
\nuc{\mub}{\mu_b}
\nuc{\muij}{\mu_{ij}}
\nuc{\mukj}{\mu_{kj}}
\nuc{\muab}{\mu_{ab}}
\nuc{\mbi}{\bar m_i}
\nuc{\mbk}{\bar m_k}
\nuc{\mbj}{\bar m_j}
\nuc{\mbij}{\bar m_{ij}}
\nuc{\mbkj}{\bar m_{kj}}

\nuc{\NoT}{N_0^T}
\nuc{\nut}{\nu_t}
\nuc{\nl}{n_l}
\nuc{\nm}{n_m}
\nuc{\Nt}{N_t}
\nuc{\NT}{N_T}
\nuc{\nT}{n_T}
\nuc{\nicj}{n_{i,j}}
\nuc{\Nu}{N_u}
\nuc{\Nucm}{N_{u,m}}
\nuc{\Nucl}{N_{u,l}}
\nuc{\Ntcl}{N_{t,l}}
\nuc{\Not}{N_0^t}
\nuc{\Notm}{N_0^{t_-}}

\nuc{\ome}{\omega}
\nuc{\Ome}{\Omega}
\nuc{\omekl}{\ome_{kl}}
\nuc{\omeal}{\ome_{al}}
\nuc{\omet}{\ome_t}
\nuc{\Omeicj}{\Ome_{i,j}}
\nuc{\Omebicj}{\bar\Ome_{i,j}}

\nuc{\pbi}{\bar p_i}
\nuc{\pbk}{\bar p_k}
\nuc{\pbj}{\bar p_j}
\nuc{\pba}{\bar p_a}
\nuc{\pbb}{\bar p_b}
\nuc{\pj}{p_j}
\nuc{\pk}{p_k}
\nuc{\pa}{p_a}
\nuc{\pb}{p_b}
\nuc{\pij}{p_{ij}}
\nuc{\pkj}{p_{kj}}
\nuc{\pab}{p_{ab}}
\nuc{\pbar}{\bar p}
\nuc{\pbij}{\bar p_{ij}}
\nuc{\pbkj}{\bar p_{kj}}
\nuc{\pbab}{\bar p_{ab}}
\nuc{\ptij}{\tilde p_{ij}}
\nuc{\pwtu}{p_{12}}
\nuc{\pbwtu}{\bar p_{12}}

\nuc{\Picj}{P_{i,j}}
\nuc{\PTicj}{P^T_{i,j}}
\nuc{\PTicipw}{P^T_{i,i+w}}
\nuc{\PTicipk}{P^T_{i,i+k}}

\nuc{\Pt}{P_t}
\nuc{\Pu}{P_u}
\nuc{\Ptcu}{P_{t,u}}
\nuc{\Ptck}{P_{t,k}}
\nuc{\PTtcu}{P^T_{t,u}}
\nuc{\PtcugT}{P_{t,u|T}}
\nuc{\Pstcu}{P^s_{t,u}}
\nuc{\Ptcugs}{P_{t,u|s}}
\nuc{\PTtctpw}{P^T_{t,t+w}}
\nuc{\PTtctpk}{P^T_{t,t+k}}
\nuc{\Ptt}{P^t_t}
\nuc{\Ptgt}{P_{t|t}}
\nuc{\PtgT}{P_{t|T}}
\nuc{\Pst}{P^s_t}
\nuc{\Ptgs}{P_{t|s}}
\nuc{\PTt}{P^T_t}
\nuc{\Ptmwt}{P^{t-1}_t}
\nuc{\Ptgtmw}{P_{t|t-1}}
\nuc{\Po}{P_0}
\nuc{\Poi}{P_0^{-1}}

\nuc{\Qbi}{\bar Q_i}
\nuc{\Qbk}{\bar Q_k}
\nuc{\Qbj}{\bar Q_j}
\nuc{\Qj}{Q_j}
\nuc{\Qij}{Q_{ij}}
\nuc{\Qkj}{Q_{kj}}
\nuc{\Qab}{Q_{ab}}
\nuc{\Qm}{Q_m}
\nuc{\Qjclcm}{Q_{j,l,m}}
\nuc{\Qbij}{\bar Q_{ij}}
\nuc{\Qwtu}{Q_{12}}
\nuc{\Qbwtu}{\bar Q_{12}}
\nuc{\qbi}{\bar q_i}
\nuc{\qbk}{\bar q_k}
\nuc{\qbj}{\bar q_j}
\nuc{\qba}{\bar q_a}
\nuc{\qbb}{\bar q_b}
\nuc{\qj}{q_j}
\nuc{\qij}{q_{ij}}
\nuc{\qbij}{\bar q_{ij}}
\nuc{\qkj}{q_{kj}}
\nuc{\qbkj}{\bar q_{kj}}
\nuc{\qab}{q_{ab}}
\nuc{\qbab}{\bar q_{ab}}
\nuc{\qwtu}{q_{12}}
\nuc{\qbwtu}{\bar q_{12}}

\nuc{\Rbi}{\bar R_i}
\nuc{\Rbj}{\bar R_j}
\nuc{\Ri}{R_i}
\nuc{\Rj}{R_j}
\nuc{\Rt}{R_t}
\nuc{\Rij}{R_{ij}}
\nuc{\Rkj}{R_{kj}}
\nuc{\Rab}{R_{ab}}
\nuc{\Rjclcm}{R_{j,l,m}}
\nuc{\rij}{r_{ij}}
\nuc{\Rbij}{\bar R_{ij}}
\nuc{\Rbkj}{\bar R_{kj}}
\nuc{\Rbab}{\bar R_{ab}}
\nuc{\Rwtu}{R_{12}}
\nuc{\Rbwtu}{\bar R_{12}}
\nuc{\rhoi}{\rho_{i}}
\nuc{\rhok}{\rho_{k}}
\nuc{\rhoj}{\rho_{j}}
\nuc{\rhoa}{\rho_{a}}
\nuc{\rhob}{\rho_{b}}
\nuc{\rhoij}{\rho_{ij}}
\nuc{\rhokj}{\rho_{kj}}
\nuc{\rhoab}{\rho_{ab}}
\nuc{\Rect}{R_{e,t}}
\nuc{\ra}{\rightarrow}

\nuc{\Si}{S_i}
\nuc{\Sk}{S_k}
\nuc{\Sj}{S_j}
\nuc{\Sa}{S_a}
\nuc{\Sb}{S_b}
\nuc{\Sbar}{\bar S}
\nuc{\Sbi}{\bar S_i}
\nuc{\Sbk}{\bar S_k}
\nuc{\Sbj}{\bar S_j}
\nuc{\Sba}{\bar S_a}
\nuc{\Sbb}{\bar S_b}
\nuc{\Sij}{S_{ij}}
\nuc{\Sbij}{\bar S_{ij}}
\nuc{\Skj}{S_{kj}}
\nuc{\Sbkj}{\bar S_{kj}}
\nuc{\Sab}{S_{ab}}
\nuc{\Sbab}{\bar S_{ab}}
\nuc{\Swtu}{S_{12}}
\nuc{\Sww}{S_{11}}
\nuc{\Swo}{S_{10}}
\nuc{\Sow}{S_{01}}
\nuc{\Soo}{S_{00}}
\nuc{\Swwf}{S_{11}^f}
\nuc{\Swof}{S_{10}^f}
\nuc{\Sowf}{S_{01}^f}
\nuc{\Soof}{S_{00}^f}
\nuc{\Swwb}{S_{11}^b}
\nuc{\Swob}{S_{10}^b}
\nuc{\Sowb}{S_{01}^b}
\nuc{\Soob}{S_{00}^b}
\nuc{\Sbwtu}{\bar S_{12}}
\nuc{\Sigxx}{\Sigma{xx}}
\nuc{\Sigx}{\Sigma{x}}
\nuc{\Sigy}{\Sigma{y}}
\nuc{\Sigxy}{\Sigma{xy}}
\nuc{\Sigyx}{\Sigma{yx}}
\nuc{\Sigyy}{\Sigma{yy}}

\nuc{\Sx}{S_{x}}
\nuc{\Sy}{S_{x}}
\nuc{\Sxx}{S_{xx}}
\nuc{\Sxy}{S_{xy}}
\nuc{\Syx}{S_{yx}}
\nuc{\Syy}{S_{yy}}
\nuc{\Szz}{S_{zz}}
\nuc{\Szx}{S_{zx}}
\nuc{\Sxxf}{S_{xx}^f}
\nuc{\Sxyf}{S_{xy}^f}
\nuc{\Syxf}{S_{yx}^f}
\nuc{\Syyf}{S_{yy}^f}
\nuc{\Sxxb}{S_{xx}^b}
\nuc{\Sxyb}{S_{xy}^b}
\nuc{\Syxb}{S_{yx}^b}
\nuc{\Syyb}{S_{yy}^b}
\nuc{\Sjclcm}{S_{j,l,m}}
\nuc{\sigi}{\sigma_i}
\nuc{\sigk}{\sigma_k}
\nuc{\siga}{\sigma_a}
\nuc{\sigb}{\sigma_b}

\nuc{\tn}{t_n}
\nuc{\tr}{t_r}
\nuc{\ti}{t_i}
\nuc{\ticl}{t_{i,l}}
\nuc{\trcm}{t_{r,m}}
\nuc{\twnm}{t_1^{n_m}}
\nuc{\twnl}{t_1^{n_l}}
\nuc{\twclnl}{t_{1,l}^{n_l}}
\nuc{\trmw}{t_{r-1}}
\nuc{\trpw}{t_{r+1}}
\nuc{\thehw}{\hat\theta_1}
\nuc{\theho}{\hat\theta_0}
\nuc{\thek}{\theta_k}
\nuc{\thej}{\theta_j}
\nuc{\thea}{\theta_a}
\nuc{\theb}{\theta_b}
\nuc{\theab}{\theta_{ab}}
\nuc{\thekj}{\theta_{kj}}
\nuc{\dtri}{\tr - \ti}
\nuc{\dtrcmicl}{\trcm - \ticl}
\nuc{\dTtr}{T - \tr}
\nuc{\dTti}{T - \ti}
\nuc{\dTticl}{T - \ticl}
\nuc{\dtr}{\trpw-\tr}
\nuc{\dtrcm}{t_{r+1,m}-\trcm}
\nuc{\tauwnicj}{\tau_1^{\nicj}}
\nuc{\tauicj}{\tau_{i,j}}
\nuc{\tauicjcl}{\tau_{i,j,l}}

\nuc{\ui}{u_i}
\nuc{\uk}{u_k}
\nuc{\uj}{u_j}
\nuc{\uij}{u_{ij}}
\nuc{\ukj}{u_{kj}}
\nuc{\uab}{u_{ab}}
\nuc{\ubi}{\bar u_i}
\nuc{\ubj}{\bar u_j}
\nuc{\Ui}{U_i}
\nuc{\Uj}{U_j}
\nuc{\Uij}{U_{ij}}
\nuc{\Ubi}{\bar U_i}
\nuc{\Ubj}{\bar U_j}

\nuc{\vi}{v_i}
\nuc{\vj}{v_j}
\nuc{\vij}{v_{ij}}
\nuc{\vbi}{\bar v_i}
\nuc{\vbj}{\bar v_j}
\nuc{\Vi}{V_i}
\nuc{\Vj}{V_j}
\nuc{\Vij}{V_{ij}}
\nuc{\Vbi}{\bar V_i}
\nuc{\Vbj}{\bar V_j}

\nuc{\wt}{w_t}
\nuc{\wtmw}{w_{t-1}}
\nuc{\wicj}{w_{i,j}}
\nuc{\wicjclcm}{w_{i,j,l,m}}

\nuc{\xt}{x_t}
\nuc{\xtpw}{x_{t+1}}
\nuc{\xtmw}{x_{t-1}}
\nuc{\xo}{x_0}
\nuc{\xw}{x_1}
\nuc{\xto}{x_2}
\nuc{\xT}{x_T}

\nuc{\xTt}{x^T_t}
\nuc{\xtt}{x^t_t}
\nuc{\xtmwt}{x^{t-1}_t}

\nuc{\xtil}{\tilde{x}}
\nuc{\xtilt}{\tilde{x}_t}
\nuc{\xtilu}{\tilde{x}_u}
\nuc{\xtiltgs}{\tilde{x}_{t|s}}
\nuc{\xtilugs}{\tilde{x}_{u|s}}
\nuc{\xtiltgt}{\tilde{x}_{t|t}}
\nuc{\xtiltgT}{\tilde{x}_{t|T}}
\nuc{\xtilugT}{\tilde{x}_{u|T}}
\nuc{\xtiltgtmw}{\tilde{x}_{t|t-1}}
\nuc{\xhat}{\hat{x}}
\nuc{\xhatt}{\hat{x}_t}
\nuc{\xhattgs}{\hat{x}_{t|s}}
\nuc{\xhattgtmw}{\hat{x}_{t|t-1}}
\nuc{\xhattt}{\hat{x}_t^t}
\nuc{\xhattgt}{\hat{x}_{t|t}}
\nuc{\xhattgT}{\hat{x}_{t|T}}
\nuc{\xhattmwgtmw}{\xhat_{t-1|t-1}}
\nuc{\xtgTf}{x_{t|T}^f}
\nuc{\xtgTb}{x_{t|T}^b}
\nuc{\xtmwgTf}{x_{t-1|T}^f}
\nuc{\xtmwgTb}{x_{t-1|T}^b}
\nuc{\ytgTf}{y_{t|T}^f}
\nuc{\ytgTb}{y_{t|T}^b}
\nuc{\ytmwgTf}{y_{t-1|T}^f}
\nuc{\ytmwgTb}{y_{t-1|T}^b}

\nuc{\yt}{y_t}
\nuc{\yr}{y_r}
\nuc{\ywn}{y_1^n}
\nuc{\ywT}{y_1^T}
\nuc{\yws}{y_1^s}

\nuc{\zwT}{z_1^{2T}}

\nuc{\ben}{\begin{enumerate}}
\nuc{\een}{\end{enumerate}}
\nuc{\itm}{\item[]}
\nuc{\bult}{$\bullet$}
\nuc{\mbo}{$\mbox{ }$}

\nuc{\calj}{{\cal J}}
\nuc{\callTtn}{\ell(T,n;t_1^n)}
\nuc{\cHT}{\cH_T}

\nuc{\calhr}{{\cal H}_R}
\nuc{\calhk}{{\cal H}_K}
\nuc{\emu}[1]{e^{-#1}}

\nuc{\LamT}{\Lam(T)}
\nuc{\LT}{L_T}
\nuc{\LTtn}{L(T,n;t_1^n)}
\nuc{\Ldk}{L_k}
\nuc{\Ldj}{L_j}
\nuc{\LkTtn}{L_k(T,n;t_1^n)}

\nuc{\Piwn}{\Pi_1^n}
\nuc{\pit}{{\it Proof}. }
\nuc{\phiT}{\phi_T}

\nuc{\pik}{\pi_k}
\nuc{\Pe}{P_e}

\nuc{\rhoT}{\rho_T}

\nuc{\sfracwt}{\frac{1}{2}}
\nuc{\sfracwf}{\frac{1}{4}}
\nuc{\sfracws}{\frac{1}{6}}
\nuc{\bha}{Bahattacharya }
\nuc{\afy}{affinitiy }
\nuc{\afY}{affinity}

\nuc{\ppS}{point processes}
\nuc{\pps}{point processes }
\nuc{\pP}{point process}

\nuc{\rps}{renewal processes }
\nuc{\rpS}{renewal processes}
\nuc{\rit}{{\it Remark}. }

\newcommand{\Bpf}{\begin{pf*}}
\newcommand{\Epf}{{\hfill$\square$}\end{pf*}}

\nuc{\fR}{\mathfrak{R}}
\nuc{\fL}{\mathfrak{L}}
\nuc{\fH}{\mathfrak{H}}
\nuc{\cLk}{\cL_k}

\nuc{\PeT}{P_{e,T}}
\nuc{\cJ}{\mathcal{J}}

\newtheorem{theorem}{Theorem}

\newtheorem{lemma}{Lemma}
\newtheorem{result}{Result}
\newtheorem{definition}{Definition}
\newtheorem*{example}{Example}

\begin{document}

\title{Asymptotic Error Rates for \\Point Process Classification
\footnote{
This work has been submitted to the IEEE for possible publication. 
Copyright may be transferred withour notice, 
after which this version may no longer be accessible.
}}

\author{Xinhui Rong
	and Victor Solo
        \thanks{School of Electrical Engineering \& Telecommunictions, UNSW, Sydney.}}

\maketitle

\begin{abstract}
Point processes are finding growing applications in numerous fields, 
such as neuroscience, high frequency finance and social media.
So classic problems of
classification and clustering are of increasing interest.
However, analytic study of  misclassification error probability
in multi-class classification has barely begun.
In this paper, we tackle the multi-class likelihood classification problem 
for point processes 
and develop, for the first time, both asymptotic upper and lower bounds on the error rate 
in terms of computable pair-wise affinities.
We apply these general results to 
classifying renewal processes. 
Under some technical conditions, 
we show that the bounds have exponential decay
and give explicit associated constants.
The results are illustrated with a non-trivial simulation.
\end{abstract}

\section{\bf Introduction}
In the past several decades,
point processes\footnote{A point process is just the sequence 
of times at which an event of interest occurs.}
 (a.k.a. spike trains) 
have arisen in a wide range of applications, such as neural coding \cite{William97}, 
genomics \cite{Sandelin10}, seismology \cite{Schoenberg08}. 
More recently, analysis of point process data has revealed the temporal dynamics 
of social media data \cite{Zhou13} and has found a potential application 
in event triggered state estimation \cite{Shi15}. 
There is also widespread interest in point fields (i.e. where the points
are distributed in space rather than along a time line)
\footnote{Also misleadingly called spatial point processes;
but if time is not involved, the word `process' does not apply.}
but we do not consider them here.

Now the emergence of large scale point process data
has led to an interest in the
classification and clustering of point process trajectories based on their statistical properties.
Lukasik et al. \cite{Lukasik16} modeled Twitter data with Hawkes processes. 
Victor and Purpura \cite{Victor97} introduced a point-wise distance between point processes 
 to cluster neural spike trains from the auditory and visual cortex.

Despite this emerging interest 
basic properties such as classification error probability (a.k.a. error rate) 
are barely discussed. 
Our previous work \cite{Rong21} was the first to tackle the problem 
but only treated binary classification (i.e. classification into one of two groups)
and only found an upper bound.
There was no asymptotic analysis.
Pawlak et al. \cite{Pawl23a, Pawl23b} studied
binary classification of time-varying Poisson processes, 
and developed error rate bounds
for various classification schemes.
Their analysis included some asymptotics.

Finding bounds for the multi-class (MC) error rate
is challenging, 
even in non point process scenarios. 
For binary classification, some upper bounds are well-known \cite{Duda01}, 
e.g. Chernoff bound \cite{Chernoff52}, Bhattacharya bound \cite{Kailath67} and 
Shannon bound \cite{Kovalevskij65}. More recent binary work
\cite{Ding22} uses
mutual information.
But for MC problems one must move to 
entropy or (Kullback-Liebler) divergence, 
e.g. \cite{Chen76, Lin91, Prasad15, Wisler16, Sekeh20}, 
however, they seldom yield analytical expressions \cite{Wisler16,Sekeh20}. 
None of these works have asymptotic analysis of error rates. 
There is asymptotic analysis in \cite{Blah74, Davi81, Gutm89, Unni16} 
but the random variables are discrete valued  (multinomial) distributions. 
Another apparently related domain is source coding where
MC classification has been discussed.
But all this deals with finite alphabets in discrete time.
We can get an approximate representation of a point process
in discrete time by binning the counts,
but the resulting alphabet is not finite, since any count, no matter how large has
a positive probability of occurring.
Finally then, none of the work in this paragraph deals with point processes.

In this paper, we tackle, for the first time, the challenging problem of finding
both asymptotic upper and lower bounds on the error rate 
for classifying point processes from MCs. We consider 
likelihood based classification (aka Bayes rule). 
We derive general bounds based on pair-wise affinities
and an asymptotic inverse Fano theorem.
We then apply these results to renewal processes.
We show that the error rate bounds have asymptotic exponential decay 
under some technical conditions. 
We give explicit asymptotic rates and amplitude constants.

Since renewal process event times are sums of 
independent and identically distributed (iid) inter-event times, 
one might suppose the classification problem reduces to iid classification,
which is of course well studied.
But there is a huge difference since iid classification deals with a fixed number of events
whereas renewal process classification deals with a fixed observation time interval
and so a variable number of events.
This makes the analysis of these two cases totally different.

The rest of the paper is organized as follows. 
In Section \ref{prel}, we review point process and likelihood classification 
properties. 
In Section \ref{bd}, we find upper and lower bounds 
on the error rate based on 
pairwise affinities.
Additionally we provide 
for the first time, an asymptotic inverse Fano theorem.
In Section \ref{aff}, we derive the asymptotics of pair-wise affinities 
for renewal processes.
In section \ref{main} 
we put all this together to get 
asymptotic error rate bounds for MC classification of renewal processes.
In Section \ref{sim}, we show comparative simulation results. 
Section \ref{con} contains conclusions. 
There are three appendices which contain proofs. 

\section{\bf Preliminaries}
\setcounter{equation}{0}
\label{prel}
We review point processes and likelihood classification properties. 
We shall see that the characterization of 
\pp es involve a hybrid likelihood of both discrete and continuous random variables 
and thus, the likelihood classification properties differ from the usual. 

\subsection{\bf Point Process Properties}
A scalar point process
can be described in 
two equivalent ways:
\nlist{
\itum{i} as a sequence of event times $T_1<T_2<\dotsm<T$. 
\itum{ii} as a counting process: 
$N_t = N_{(0,t]} = $no. of events up to and including time $t$. 
}
A point process is characterized by its 
intensity function, which is the probability of 
getting a new event in the next small time interval 
given history up to time $t$ and has the form
\eq{
\lam(t) = \lam(t|\Notm)= \lim_{\delta\to0} \frac1\delta P[N_{t+\delta} - N_t=1|\Notm],
}
where $\Notm= \{N_u, 0<u< t\}$ is the history up to current time $t$. 
Generally, the intensity is stochastic; but for Poisson processes it is deterministic.

Renewal processes (RPs) are characterized by having 
iid inter-event times (IETs) $X_i=T_i-T_{i-1}$
with common IET density $p(x)$. 
The next event time is dependent on the immediately preceding event time. 
Introduce the cumulative distribution function (cdf) $F(x)=\int_0^x p(u)du$,
the survivor function $S(x)=1-F(x)$, hazard function 
$h(x)=\frac{p(x)}{S(x)}$ and integrated hazard $H(x) = \int_0^x h(u)du$.
Then the stochastic intensity 
 is given by \cite{Daley03}
\eq{
\lam(t) = h(t-T_r), \quad T_r\leq t<T_{r+1}.
}
The classic hazard relations \cite{Daley03}
\eq{
p(x) = h(x)e^{-H(x)}, \quad S(x) = e^{-H(x)}
}
are frequently used in later sections. 

We assume throughout that $0<h(x)<\infty$. This avoids
pathological cases. 
Also introduce the integrated intensity $\Lam(t) = \intot \lam(u)du$. 
It is straightforward to show that on $[T_r, T_{r+1})$
\eq{
\arr{lcl}{
\Lam(t) &=& \Lam(T_r) + H(t-T_r) \mbox{ where}\\
\Lam(T_r) &=& \ssum{k=1}^r H(X_k).}
}

\noindent\textbf{(iii)} The likelihood / \textit{Janossy density}:

The likelihood function of a \pp\ is a hybrid density of 
the total number of events $N_T$, and event times $T_1^{N_T}$, 
observed in time period $[0,T]$. 

If we denote a trajectory $\NoT=\{N_T=n, T_1^{N_T}=\twn=\{t_1,\dotsm,t_n\}\}$, 
then, the likelihood (aka the \textit{Janossy density}) has the form \cite{Daley03}
\eq{
\arr{rcll}{
&L(\NoT) &=& \LTtn \\
\mbox{if }\NT=0&&= &
P(N_T=0) = e^{-\Lambda(T)}\\
\mbox{if }\NT>0&&=&
P(N_T=n,T_1^n=t_1^n)\\
&&=&
e^{-\LamT}\prod_{r=1}^n\lam(t_r)
}}
The Janossy density sums, integrates to $1$ as follows \cite{Daley03},
\eq{
1 &= P[N_T=0] + \ssum{n=1}^\infty P[N_T=n]\\
	&= e^{-\Lambda(T)} +\ssum{n=1}^\infty\int_{R_n(T)} L(T,n;\twn)d\twn\\
R_n(T) &= \{\twn|0<t_1<t_2<\dotsm<t_n<T\}.
}
For simplicity, in the rest of the paper, we informally write the sum/integral 
$\ints \LTtn d\twn=\ints L(\NoT) d\twn$.

\subsection{\bf Likelihood Classification}
\subsubsection{Likelihood Rule/Bayes Rule}

Consider the problem of classifying a \pp\ trajectory into one of $M\geq 2$ classes. 
Class $k$ has intensity function $\lam_k(t)$ and occurs with prior probability $\pik$.
The class label $C$ is a random variable with mass function,
$P(C=k)=\pik$ and $\ssum{k=1}^M \pik=1$.

Then the class $k$ log-likelihoods are
\eq{
\ell_k(\NoT) = \ln L_k(\NoT)=
\cas{
-\Lamk(T),&N_T=0\\
\ln\lamk(\tr) -\Lamk(T),&N_T>0.
}}

The likelihood rule assigns the trajectory $\NoT$ to class $k$ if
\eq{
\ln\pi_k + \ell_k(\NoT) \geq \ln\pi_j + \ell_j(\NoT), \mbox{for all }j\neq k.
}

We also use the mixture likelihood
\eq{
L=\ssum{k=1}^M \pik L_k(\NoT),
}
and the class expectation notation
\eq{
\E_k[f(N_T,T_1^{N_T})] &= \ints f(n,\twn)L_k(\NoT) d\twn.}

\subsubsection{Misclassification Error Rate}

Firstly we consider the classification problem 
for analog data $X$. 
A classifier, i.e. an estimator of the class label $\hat C(X)$, 
has misclassification error \cite{Devroye96}
\eq{
e = \skwM \pi_k P[\hat C(X)\neq k|C=k]. 
}
This is minimized by the likelihood classification rule, 
yielding misclassification error probability
\eq{
P_e =  1 - \int \max\{\pi_1 L_1(x),\dotsm,\pi_M L_M(x)\}dx.
}
In the \pp\ case, we have Janossy densities, thus
\eq{
P_e(T) =  1 - \ints \max\{\pi_1 L_1,\dotsm,\pi_M L_M\}d\twn,
}
which we note, is a function of observation time $T$.

\section{\bf Error Rate Bounds: Affinities and Asymptotics}
\setcounter{equation}{0}
\label{bd}
In this section, we consider calculating the mis-classification error rate.
In the first subsection we develop new upper and lower error rate bounds
based on Shannon entropy.
But entropies turn out to be very hard to calculate.
So in the second subsection,
we extend to point processes, the recent method of \cite{Kolchinsky17}
which  bounds the entropies in terms of pair-wise affinities.
But the affinities can only be calculated asymptotically.
So  in the third
subsection, we develop a general asymptotic
 inverse Fano theorem which will enable
 asymptotic evaluation of the affinity bounds.

\subsection{Point Process Error Rate Entropy Bounds}
We first define the Shannon entropy 
and mutual information as applied to \pp\ classification as follows. 

\defi{
The Kullback-Liebler (KL) divergence $D(f||g)$ between 
two \pp\ densities $f,g$ has the form
\eq{
D(f||g) = \ints f(\NoT)\ln\frac{f(\NoT)}{g(\NoT)}d\twn.
}
}

\defi{
(a) The Shannon entropy is
\eq{
\cH(T)&= \cH(C|\NoT) \\
	&= -\E\sqbra{\skwM P[C=k|\NoT]\log P[C=k|\NoT]},
}
where $\log$ is $\log_2$ and $\E[\cdot]$ is expectation taken regarding the mixture, 
throughout the paper. 

(b) The mutual information (MI) is 
\eq{
\cI(T) =\cI(C;\NoT) &= \skwM\pi_k \E_k\sqbra{\log\frac{L_k}{L}}\\
	&=\skwM \frac{\pi_k}{\ln2}D(L_k||L).
}
}

They are both functions of observation time $T$. 
Also note that in terms of entropies the MI is
\eq{
\cI(C;\NoT) = \cH(C) - \cH(C|\NoT),
}
where $\cH(C) = -\skwM\pi_k\log\pi_k$ is the discrete entropy of priors. 

We now have the following result.

\lem{
\label{ulb1}
{Multiclass Error Rate Entropy Bounds.}
The misclassification error rate $P_e$ is bounded as follows.
\eq{
\phi(\cH(T))\leq P_e(T)\leq  \wo2\cH(T),
}
where $\phiT=\phi(\cH(T))$ solves
\eq{
&\cH_b(\phiT) +  \phiT\log(M-1) = \cH(T),
}
and $\cH_b(\phi) = -\phi\log \phi - (1-\phi) \log(1-\phi)$ is the binary entropy.
}
{\it Proof}.
For analog random variables $X\in\mathbb R^n$, 
the proof of the upper bound can be found in \cite{Kovalevskij65},\cite{Lin91} 
and the lower bound is a direct application 
of Fano's inequality \cite{Cover91}. 
The proofs both extend to \pp es without difficulty. \hfill$\square$

{\it Remarks}.

(i) Some other upper bounds given in \cite{Kovalevskij65} might be tighter when $T$ is small. 
But for large $T$ the
bound quoted above is the tightest. 

(ii) However, the Shannon entropy $\cH(T)$, generally has no analytic expression and 
finding its asymptotics directly is also notoriously difficult 
given the unusual likelihood of \pp es. 
So we follow \cite{Kolchinsky17} and bound it
in terms of pairwise affinities by
developing a point process version of \cite{Kolchinsky17}.

\subsection{Point Process Error Rate Affinity Bounds}
We first introduce the following affinities for \pp es. 
\defi{
The Bhattacharyya affinity has form
\eq{
\Rkj(T) &= \ints\sqrt{L_k(\NoT)L_j(\NoT)}d\twn, \\
\intertext{and the KL affinity has form}
\Kkj(T) &=e^{-\Dkj(T)}, \quad \Dkj(T) = D(L_k||L_j).
}
}
Note that $-\ln\Rkj$ and $-\ln\Kkj$ are the \Bhatt\ divergence 
and the KL divergence, respectively.
We give the \pp\ bounds in the following theorem.

\lem{
\label{ulb2}
Entropy and MI Bounds via Afinities.

(a) The MI is bounded as follows
\eq{
M_R(T) \leq \cI(T) \leq M_K(T),
}
where $M_R(T) = -\skwM\pi_k\log\sjwM\pi_j\Rkj(T)$ and \\
$M_K(T) = -\skwM\pi_k\log\sjwM\pi_j\Kkj(T)$.

(b) The Shannon entropy is bounded as follows
\eq{
\cH_K(T)\leq\cH(T)\leq\cH_R(T),
}
where 
\eq{
\cH_R(T) &= \skwM\pi_k\log\sqbra{1+\ssum{j\neq k} \frac{\pi_j}{\pi_k} \Rkj(T)}\\
\cH_K(T) &= \skwM\pi_k\log\sqbra{1+\ssum{j\neq k} \frac{\pi_j}{\pi_k} \Kkj(T)}.
}

(c) The misclassification error probability $P_e(T)$ is then bounded as follows
\eqn{
\label{bds}
\phi(\cH_K(T))\leq P_e(T)\leq \wo2\cH_R(T).
}
}

{\it Proof}. See appendix A.

\subsection{Asymptotic Inverse Fano Theorem}
We first define the asymptotic decay rate.
We say a function $f(T)$ asymptotes to $G(T)$ at $\tau$, 
or $f(T) \sim G(T)$ at $\tau$, with $\tau$ to be finite or infinite, if
$\lim_{T\ra \tau} \frac{f(T)}{G(T)} = 1.    $

To develop the asymptotic bounds related to the affinities, 
we need the following asymptotic inverse Fano theorem. 

\theo{
\label{invF}
Asymptotic Inverse Fano Theorem. \\
Suppose $0<\phiT<1$ solves
\eq{
-\phiT\log\phiT - (1-\phiT)\log(1-\phiT) + d \phiT = \cH(T),
}
where $d$ is a constant and $\cH(T)\downarrow0$, as $\Ttoi$. Then,
\eq{
\phiT\sim\frac{\cH(T)}{-\log\cH(T)}\mbox{ as } T\toi.
}
}

Theorem \ref{invF} enables asymptotics of the lower bound. 
Then using the equivalence $\log(1+f)\sim\frac1{\ln2}f$ as $f\to0$, 
it is straightforward to find that the bounds in (\ref{bds}) 
have the following asymptotic behaviors 
related to the affinities
\eq{
\arr{rrcl}{
	&(\ln2)\calhr(T)&\sim& \ssum{k=1}^M\pik\ssum{j\neq k}\pi_j\Rkj(T)\\
	&(\ln2)\calhk(T)&\sim& \ssum{k=1}^M\pik\ssum{j\neq k}\pi_j\Kkj(T)\\
\Ra &\phiT=\phi(\calhk(T))&\sim& \frac{\calhk(T)}{-\log\calhk(T)}, \mbox{ as }\Ttoi.
}}

\section{\bf Asymptotic Affinities for Renewal Processes}
\setcounter{equation}{0}
\label{aff}
Here, we consider classification for RPs and 
develop the asymptotic affinities based on Laplace transform (LT) analysis. 
In the first subsection derive the LTs of the affinities and divergences. 
Then in the following two subsections, 
we derive their asymptotics.

\subsection{Laplace Transform of Affinities and Divergences}
We first derive the various LTs of 
the \Bhatt\ affinity $\Rkj$ and KL divergence $\Dkj$ for $k\neq j$. 
We denote $\bar f(s)=\intoi f(t)e^{-st}dt$ the LT of $f(t)$. 

Recall from Section \ref{prel} A that 
for class $k$, the IET has density $p_k(x)$, hazard $\hk(x)$, survivor $\Sk(x)$ 
and integrated hazard $\Hk(x)$. 
We introduce
\eq{
\arr{rclll}{
\pkj(x) &=& \sqrt{\pk(x)\pj(x)}, \quad \Gkj(x) &=& \sqrt{\Sk(x)\Sj(x)}\\
 \qkj(x) &=& \pk(x)\ln\frac{\pk(x)}{\pj(x)}, \quad \Fkj(x) &=& \Sk(x) \ln \frac{\Sk(x)}{\Sj(x)}.
}}
We have the following LT results. 

\lem{  
\label{RPub}
The  \Bhatt\ affinity $\Rkj(T)$ has LT
\eq{
\Rbkj(s) = \frac{\bar{G}_{kj}(s)}{1-\pbkj(s)}.
}
}

{\it Proof}.
See \cite{Rong21}.

\lem{ 
\label{RPlb}
The KL divergence 
$\Dkj(T)$ has LT
\eq{
\Dbkj(s) &= \frac{\Bbkj(s)}{s(1-\pbk(s))} = \frac{\qbkj(s) + s\Fbkj(s)}{s(1-\pbk(s))}.
}
}

{\it Proof}. See appendix B.

{\it Remark}.

To recover time domain expressions, 
numerical inversion is possible in some cases, 
e.g. \cite{Gaver66}\cite{Stehfest70}\cite{Dubner68}\cite{Talbot79}\cite{Horvath20}. 
However, we have shown in \cite{Rong21} that all methods have problems 
when the bounds are close to $0$, 
i.e. when observation time $T$ is large. 
So we proceed to derive the analytical asymptotics. 

\subsection{$\Rkj$ Asymptotics}
We assume regular IET densities 
whose moment generating functions (MGFs) exist as defined below. 
Heavy-tailed distributions do not have MGFs and will be treated elsewhere.

\defi{
\label{mgf}
We call a density $p$ regular, iff there exists a real number 
$\theta^\ast>0$,  such that
\eq{
m(\theta) = \pbar(-\theta) = \intoi e^{\theta x} p(x)dx < \infty,
}
for all $\theta<\theta^\ast$, and $m(\infty)=\infty$.
}

We now obtain an exponential decay formula for 
the \Bhatt\ affinity $\Rkj(T)$ under two conditions.

\theo{
\label{Ra}
$\Rkj$ Asymptotics. 
Suppose $\pk,\pj$ are both regular. Introduce $\pkj^*(x)=e^{\gamkj x}\pkj(x)$. 
Then, we can find a unique $\gamkj>0$ to make $\pkj^*(x)$ a density i.e.
$\intoi \pkj^*(x)dx =1$.

If (a) $\Gbkj(-\gamkj) =\intoi e^{\gamkj x}\Gkj(x)dx <\infty$ and 
(b) $-\pbkj'(-\gamkj) =\intoi x \pkj^*(x)dx<\infty$, 
then as $T\toi$,
\eqn{
\label{Rkja}
\Rkj(T) \sim \alpkj e^{-\gamkj T}, 
}
where
\eq{
\alpkj = \frac{\Gbkj(-\gamkj)}{-\pbkj'(-\gamkj)} 
=
\frac{\intoi e^{\gamkj x} \Gkj(x) dx}{\intoi x \pkj^*(x)dx}.
}}

{\it Proof}. See appendix B.

Now we are left with two conditions to be checked:
\nlist{
\ita $\intoi e^{\gamkj x}\Gkj(x)dx<\infty$ and, 
\itb $\intoi x \pkj^*(x)dx<\infty$
}

We show in Appendix B that condition (b) is almost surely satisfied 
for regular densities. 
We now give a general sufficient condition under which (a) is satisfied. 

{\bf Assumption A1. } 
Hazard bounded away from $0$. 
\eq{
\liminf_{x\toi} h(x)>0.
}

A1 is a general condition that almost all frequently used regular distributions obey, 
e.g. (mixture of) exponential, gamma, normal distributions. 
We now show condition (a) is satisfied under A1 in the following lemma.

\lem{
\label{Gp}
Suppose $\pk,\pj$ are regular and their hazards satisfy A1. 
Then, if $\pbkj(-\gamkj) = 1$, we have
\eq{
\Gbkj(-\gamkj) = \intoi e^{\gamkj x} \Gkj(x)dx <\infty.
}
}

{\it Proof}. See appenix B.

\subsection{$\Kkj$ Asymptotics}
$\Kkj(T)$ has exponential decay under the following assumptions.

{\bf A2. } 
Finite IET KL divergence. 
We assume \\
$\qbkj(0)=\intoi\qkj(x)dx=\intoi \pk(x)\ln\frac{\pk(x)}{\pj(x)}dx<\infty$.

{\bf A2*. } 
$-\qbkj'(0)=\intoi x\qkj(x)dx<\infty$.

\theo{ 
\label{Ka}
Under A2 and A2*, if both $\pk,\pj$ are regular, 
then, as $T\to\infty$, 
\eqn{
\label{Kkja}
\Kkj(T) &\sim \ckj e^{-\rhokj T}, \quad \rhokj = \frac{\qbkj(0)}{\muk},
}
where 
\eq{
\muk &=\intoi x\pk(x)dx = \intoi \Sk(x)dx=\Sbk(0) \\
\ckj &= \exp\bra{\wo{\muk}\sqbra{\qbkj'(0) + \frac12\rhokj(\sigk^2 + \muk^2) + \Fbkj(0)}}\\ 
\sigk^2 &= \var_k[X]=\intoi x^2 \pk(x)dx - \muk^2.
}
}

{\it Proof}. See appendix B.

Since we only have the LT of the KL divergence $\Dkj(T)=-\ln\Kkj(T)$, 
we use a different proof and need a stronger asymptotic relation.

\section{\bf Error Rate Asymptotics for Renewal Processes}
\setcounter{equation}{0}
\label{main}
Here we use the RP affinity asymptotics from section IV
in the general asymptotic results of section III
to get our main theorem.
We illustrate the theorem with a non-trivial
example having Gamma distributed IETs.

\theo{
\label{ab}
Error Rate Asymptotics for the Renewal Process.  
Assume A1, A2, A2$^\ast$ and regular $\pk(x)$.
Then $\Rkj(T)\sim\alpkj e^{-\gamkj T}$ and
$\Kkj(T)\sim\ckj e^{-\rhokj T}$ as given in (\ref{Rkja}) and (\ref{Kkja}), respectively.
Further, the upper and lower error bounds behave as follows:
$\phiT=\phi(\calhk(T))\leq \Pe(T)\leq \frac12\calhr(T)$ where
\eq{
\cH_R(T) &\sim \frac{\alpha^*}{\ln2} e^{-\gam^* T}\mbox{ and }
\phiT\sim \frac{c^*}{\rho^*T} e^{-\rho^* T}.
}
where $\gam^*=\min_{k\neq j}\{\gamkj\}$, $\rho^* = \min_{k\neq j}\{\rhokj\}$, \\
$\alpha^*=\ssum{\gamkj=\gam^*,k\neq j} \pi_j \alpkj$ and 
$c^*=\ssum{\rhokj=\rho^*,k\neq j} \pi_j \ckj$.
}

{\it Proof}. See appendix C.

{\bf Gamma Example. } 
Here, we apply the results to the Gamma distributed IET.
We assume class densities
\eq{
p_k(x) = p(x;\thek,\betk) = \frac{\betk^{\thek}}{\Gam(\thek)}x^{\thek-1}e^{-\betk x},
}
where $\Gam(\theta) = \intoi x^{\theta-1}e^{-x}dx$ is the gamma function.

{\it (a) Upper Bound Asymptotics }

We need to find $\gamkj,  -\pbkj'(-\gamkj),  \Gbkj(-\gamkj)$.

First, let $\thekj = \frac12(\thek+\thej), \betkj = \frac12(\betk+\betj)$ and $\dkj = \sqrt{\frac{\betk^{\thek}\betj^{\thej}}{\Gam(\thek)\Gam(\thej)}}$.
Then we can find that
\eq{
\intoi e^{\gamkj x}\pkj(x) dx &= \dkj \intoi  x^{\thekj-1}e^{-(\betkj-\gamkj) x}dx\\
	&= \dkj \frac{\Gam(\thekj)}{(\betkj-\gamkj)^{\thekj}}=1\\
\Ra \gamkj &= \betkj - (\dkj\Gam(\thekj))^{1/\thekj}.
}

Further, 
\eq{
-\pbkj'(-\gamkj) &= \intoi x e^{\gamkj x} \pkj(x)dx \\
	&= \dkj\intoi x^{\thekj}e^{-(\betkj-\gamkj)x}dx\\
	&= \dkj\frac{\Gam(\thekj+1)}{(\betkj-\gamkj)^{\thekj+1}} = \frac{\thekj}{\betkj-\gamkj}.
}

And $\Gbkj(-\gamkj) =\intoi e^{\gamkj x}\Gkj(x)dx$ 
can then be found numerically. \\

{\it (b) Lower bound asymptotics }

We need to find: $\qbkj(0),\qbkj'(0),\Fbkj(0)$
and $\E_k[X^2]=\muk^2+\sigk^2$.

It is well-known that $\muk = \frac{\thek}{\betk}, \sigk^2 = \frac{\thek}{\betk^2}$ and 
$\qbkj(0) =$
\eq{
(\thek-\thej)\psi(\thek)- \ln\frac{\Gam(\thek)}{\Gam(\thej)} + \thej\ln\frac{\betk}{\betj} + \thek\frac{\betj-\betk}{\betk},
}
where $\psi(x)=\frac{\Gam'(x)}{\Gam(x)}$ is the digamma function.

We now find
\eq{
-\qbkj'(0) &= \intoi x\pk(x)\ln\frac{\pk(x)}{\pj(x)}dx\\
	&= \intoi x\pk(x)\ln\frac{x\pk(x)}{\pj(x)}dx - \intoi x\pk(x)\ln xdx\\
	&= \muk \intoi p(x;\thek+1,\betk) \ln\frac{p(x;\thek+1,\betk)}{\pj(x)}dx\\%
	&\quad - \muk \intoi p(x;\thek+1,\betk)\ln\frac{p(x;\thek+1,\betk)}{\pk(x)}dx\\ %
		&= \muk\bra{\qbkj(0) - \frac{\thej}{\thek} + \frac{\betj}{\betk}}.
}

And $\Fbkj(0)$ can also be found numerically.

\section{\bf Simulations}
\setcounter{equation}{0}
\label{sim}

Here we show non-trivial classification simulations for RPs with 
mixture of Erlang (MoE) IET distributions and compare them with 
the asymptotic bounds.
MoE is widely used for modeling RPs \cite{Xiao19}.

We consider an $M=10$-class point process classification problem 
with IET densities
\eq{
\pk(x) = \ssum{l=1}^{m_k} \omekl \frac{\betkl ^l}{(l-1)!} x^{l-1}e^{-\betkl x}, \quad k = 1,\dotsm,10
}
where $m_k$ is the model order and $0<\omekl\leq1$ are mixture weights with 
$\ssum{l=1}^{m_k}\omekl=1$. We assume equal priors $\pi_k=1/10, k=1,\dotsm,10$.

\begin{figure}[!t]
\centering
\includegraphics[width=3.5in]{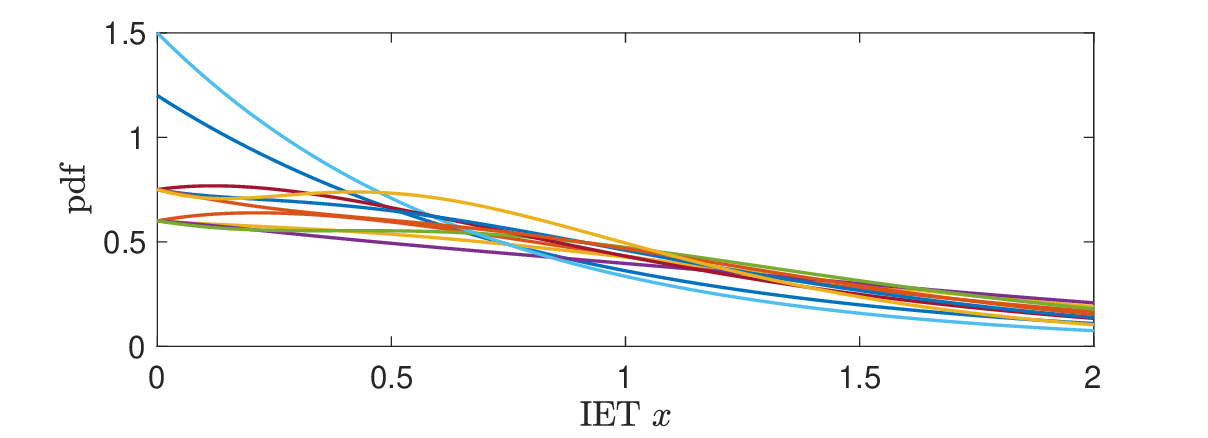}
\caption{MoE IET densities.}
\label{pdf}
\end{figure}

We specify the model parameters
\eq{
m_1^{10} &= \bmat{m_1 & \dotsm & m_{10}}^\top \\
	&= \bmat{1 & 2  & 3 & 3 & 4 & 1 & 2  & 3 & 3 & 4}^\top\\
[\omekl] &= \bmat{1 & \sfracwt & \sfracwt & \sfracwt & \sfracwt & 1 & \sfracwt & \sfracwt & \sfracwt & \sfracwt\\
	0 & \sfracwt & \sfracwf & \sfracwf & \sfracws & 0 & \sfracwt & \sfracwf & \sfracwf & \sfracws \\
	0 & 0 & \sfracwf & 0 & \sfracws & 0 & 0 & \sfracwf & 0 & \sfracws \\
	0 & 0 & 0 & \sfracwf & \sfracws & 0 & 0 & 0 & \sfracwf & \sfracws}^\top\\
[\betkl]
&= \sqbra{\smat{1.2 & 1.2 & 1.2 & 1.2 & 1.2 & 1.5 & 1.5 & 1.5 & 1.5 & 1.5\\
			- & 1.5 & 1.5 & 1.5 & 1.5 & - & 1.7 & 1.7 & 1.7 & 1.7\\
			- & - & 1.8 & - & 2 & - & - & 2.5 & - & 3\\
			- & - & - & 2 & 3.2 & - & - & - & 3 & 4.5}}^\top 
}
The maximum mixture order is $4$ and the parameters are 
chosen so that the class IET densities are close to each other. 
The IET densities are plotted in Fig. \ref{pdf}. 

Aside from $\mu_k,\sigma^2_k$, 
the parameters in asymptotic error bounds must be found numerically.
$\gamkj$ were found by varying its values so that
$\intoi x e^{\gamkj x}\pkj(x)dx$ fell within a tiny interval around $1$. 
The other parameters were found using the \texttt{vpaintegral()} command in 
MATLAB with machine arithmetic. 

For the upper bound, we find $R_{89}(T)=R_{98}(T)$ to dominate.
Then, the upper bound is
\eq{
\wo2\cH_R(T)\sim \frac{\pi_8+\pi_9}{2\ln2} \alpha_{89}e^{-\gam_{89} T} = \alpha^*e^{-\gam^*T},
}
where $\alpha^*= 0.14$ and $\gam^* = 4.8\times10^{-4}$.

For the lower bound, we find $K_{98}(T)$ to  dominate and the lower bound
\eqn{
\label{lbeq}
\phi(T)\sim\pi_8\frac{c_{98}}{\rho_{98}}\wo{T} e^{-\rho_{98} T} = \frac{c^*}{\rho^*}\wo T e^{-\rho^*T},
}
where $\frac{c^*}{\rho^*} = 76.5$ and $\rho^*= 1.9\times10^{-3}$.

To get the error probability $P_e(T)$, we run Monte Carlo (MC) simulations 
on $5$ grids between $T=5000$ and $T=10000$. 
The number of RPs simulated at each grid of $T$ is determined ${60}/{\cH_R(T)}$, 
i.e. the larger observation time $T$, the larger number of RPs.

Two RP simulation methods are possible.
The first is to use the classic inversion method \cite{Devroye86} to simulate 
IETs until the sum of the IETs exceeds $T$. 
The second is to use the thinning algorithm \cite{Lewis79}\cite{Ogata81}, 
commonly used to simulate point processes with finite intensities. 
We use the later since the first requires numerical solution for an inverse distribution.

In the `thinning' simulation
calculation of the hazard is time consuming. 
We thus pre-calculate the hazard values on a dense grid for each class 
and look up the closest values in each simulation step.
We use a similar method when calculating the log-likelihoods in the likelihood classification.

We plot the error probability $P_e(T)$, the upper bound $\wo2\cH_R(T)$ and 
the lower bound $\phi(T)$ in Fig. 2(a). 
We plot in Fig. 2(b) the transformed error $\ln P_e(T)/T$ 
and also the transformed bounds. They should converge to a constant. 
We need to point out that we plot Fig. 2(b) to show clearly that 
these quantities have asymptotically exponential decay, 
however this is a weaker version of out theoretical results, 
since our asymptotic results imply the log-asymptotics but not the inverse. 

We observe in both figures that the MC error probabilities are 
bounded by the asymptotic error bounds. 
In Fig. 2(b), we find that the error probability $P_e$ also has asymptotically 
exponential decay and the decaying rate converges closer to the upper bound rate.

\begin{figure}[t]
\label{err}
\begin{minipage}[b]{.48\linewidth}
\label{erra}
  \centering
  \centerline{\includegraphics[width=4cm]{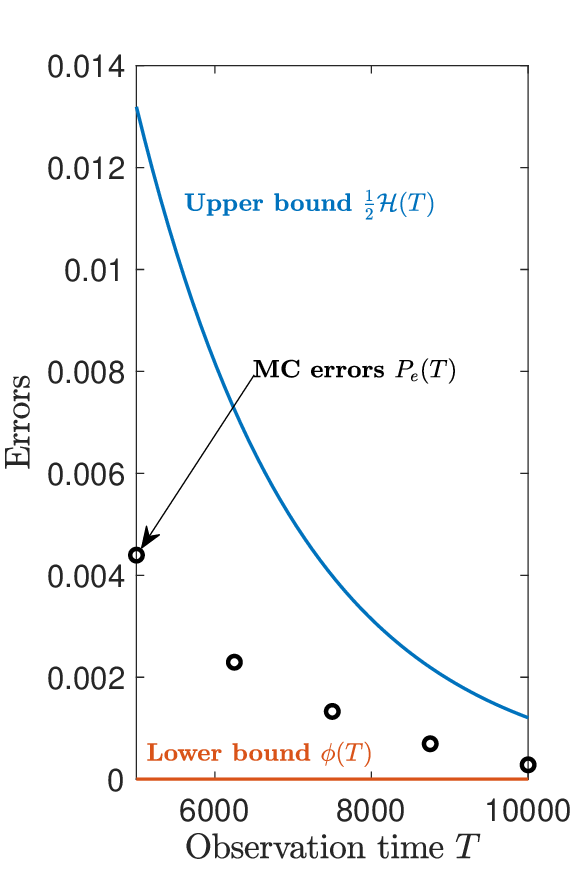}}
  \centerline{(a) Errors and bounds}\medskip
\end{minipage}
\begin{minipage}[b]{.48\linewidth}
\label{errb}
  \centering
  \centerline{\includegraphics[width=4cm]{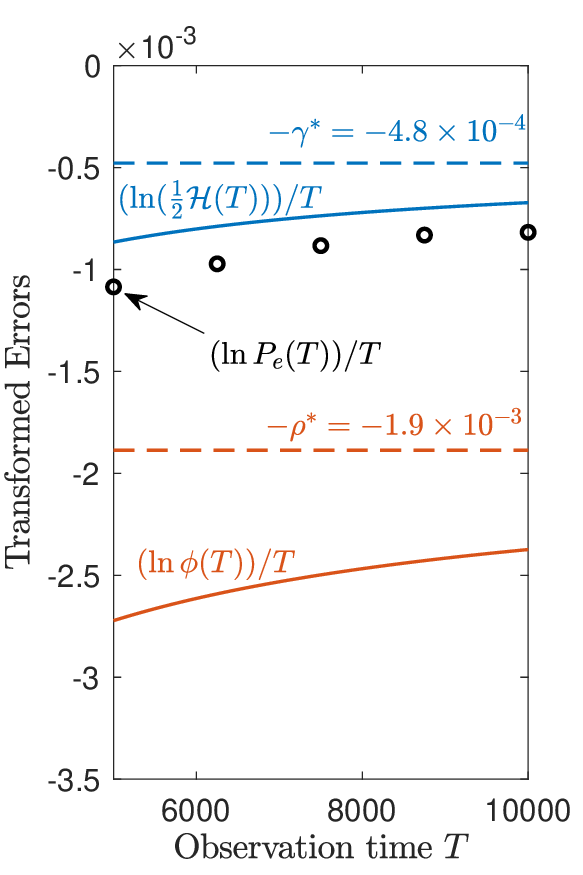}}
  \centerline{(b) {Transformed errors}{and bounds}}\medskip
\end{minipage}
\caption{MC error probabilities and error bounds: (a) For $T\geq5000$, 
the error probabilities are bounded by the asymptotic bounds. 
(b) The transformed plots show that the error probability $P_e(T)$ 
has asymptotically exponential decay.}
\label{fig:res}
\end{figure}


\section{Conclusion}
\setcounter{equation}{0}
\label{con}

In this paper we studied, for the first time,
the asymptotics of the multi-class misclassification
error probability ($P_e(T))$ for point processes.
After extending some standard entropy based analog bounds to point processes
the analysis proceeded in four stages.

Firstly we extended recent pair-wise Bhattacharyya and Kullback-Liebler
affinity based entropy bounds to the point process case.
Secondly we develop a new asymptotic inverse Fano theorem
that enables asymptotic calculation of both lower and upper bounds.
Then thirdly, we specialised to the renewal processes case.
We found the Laplace transforms of the \Bhatt\ affinity and the Kullback-Liebler divergence 
and thus derived asymptotic expressions for these affinities.

Then in the fourth stage we put  all this together
 to derive asymptotic upper and lower bounds 
on the mis-classification error probability (a.k.a. error rate)
for multi-class classification of renewal processes.
We illustrated the theory with a non-trivial Gamma example.
And further illustrated the results with a medium sized simulation using
a mixture of Erlang distributions.

In future work, we will tackle asymptotics for heavy-tailed 
renewal processes and for Hawkes processes.

\section*{Appendix A. Proofs for Section \ref{bd}}
\label{prosec3}
Here we prove Lemma \ref{ulb2} and Theorem \ref{invF}. 
The proof of Lemma \ref{ulb2} relies on the following `opposite'
Jensen's inequality. 

{\it Lemma 6}.
\label{KLub}
KL upper bound. \\
Consider the mixture density $p(x)=\skwM\pi_kp_k(x)$ 
for an analog random variable and 
introduce 
$\Dkj = D(\pk||\pj)$, $\Kkj = e^{-\Dkj}$ as well as 
$D_k = \sjwM\pi_j \Kkj$.
Then,
\eq{
D(\pk||p) \leq-\ln D_k,
}
The proof is given after the proof of Lemma 2.

\noindent
\bult\mbo{\it Proof of Lemma \ref{ulb2}.}  

(a) The MI lower bound is proved in \cite{Kolchinsky17},\cite{Ding22}. 
The MI upper bound is stated in \cite{Kolchinsky17},\cite{Ding22}, but depends on a key lemma 
from \cite{Hershey07}, \cite{Paisley10}. But \cite{Paisley10} is not available and 
\cite{Hershey07} does not have the key lemma. 
Thus, we give a new short direct proof of the upper bound 
using Lemma 6.

To finish the proof of (a) multiply through the inequality in Lemma 6
 by $\pi_k$  and sum to obtain
\eq{
\cI(T) &= \skwM \frac{\pi_k}{\ln2} D(L_k||L) \\
	&\leq -\skwM \pi_k \log(\sjwM \pi_j \Kkj) = M_K.
}

(b) 
Use the information relation 
$\cH(T) = \cH(C) - \cI(T)$ to get
\eq{
\cH(T) &\leq \cH(C) - M_R(T) \\
	&= -\skwM\pi_k\log\pi_k + \skwM\pi_k\log(\sjwM\pi_j\Rkj)\\
	&= \skwM\pi_k \log\sqbra{1+\ssum{j\neq k} \frac{\pi_j}{\pi_k}\Rkj(T)}=\cH_R(T)\\
\cH(T) &\geq \cH(C) - M_R(T) \\
	&= -\skwM\pi_k\log\pi_k + \skwM\pi_k\log(\sjwM\pi_j\Kkj)\\
	&= \skwM\pi_k \log\sqbra{1+\ssum{j\neq k} \frac{\pi_j}{\pi_k}\Kkj(T)}=\cH_K(T).
}

(c) Applying Lemma \ref{ulb1}  gives the result.\hfill$\square$\\

\bult \mbo{\it Proof of Lemma 6}. 

Let $\ome_j, j=1,\dotsm,M$ be a mass function to be chosen. 
Now with $k$ fixed, apply the log-sum inequality to find
\eq{
	\ssum{j}\bra{\ome_j p_k(x)\ln\frac{\ome_j p_k(x)}{\pi_j\pj(x)}}
\geq &\bra{\ssum{j}\ome_j\pk(x)}\ln\frac{\ssum{j}\ome_j\pk(x)}{\ssum{j}\pi_j\pj(x)}\\
= & \pk(x)\ln\frac{\pk(x)}{p(x)}.
}
Integrate both sides to obtain
\eq{
LHS &= \ssum{j} \ome_j\int\pk(x)\ln\frac{\ome_j\pk(x)}{\pi_j\pj(x)} dx 
	\geq D(\pk||p).\\
\intertext{We now rewrite the LHS as follows}
LHS &= \ssum{j}\ome_j\ln\frac{\ome_j}{\pi_j} + \ssum{j} \ome_j\int\pk(x)\ln\frac{\pk(x)}{\pj(x)}dx\\
	&= \ssum{j} \ome_j\ln\frac{\ome_j}{\pi_j} + \ssum{j} \ome_j\Dkj.
}
We choose $\ome_j$ to minimize the LHS. Adding a Lagrangian penalty 
$\lam(\ssum{j}\ome_j -1)$ and differentiating w.r.t $\ome_j$ gives
\eq{
\arr{rrcl}{
&0 &=& \ln\frac{\ome_j}{\pi_j} + 1 + \Dkj + \lam\\
\Ra &\ln\frac{\ome_j}{\pi_j} &=& -(\lam+1) - \Dkj\\
\Ra &\ome_j &=& C\pi_j e^{-\Dkj} = C\pi_j\Kkj,\quad C = e^{-(\lam+1)}.
}}
Applying the constraint $\ssum{j}\ome_j=1$ gives
$\ome_j= \frac{\pi_j \Kkj}{D_k}.$

The minimized LHS is then
\eq{
LHS &= \ssum{j} \ome_j\ln\frac{\ome_j}{\pi_j} + \ssum{j} \ome_j\Dkj\\
	&= \wo{D_k} \ssum{j} \pi_j \Kkj\ln\frac{e^{-\Dkj}}{D_k} + \wo{D_k}\ssum{j}\Dkj \pi_j\Kkj\\
	&= \wo{D_k} \ssum{j} \pi_j \Kkj\bra{-\Dkj-\ln{D_k} }+ \wo{D_k}\ssum{j}\Dkj \pi_j\Kkj\\
	&= -\ln D_k.
}
which gives the result.

\bult\mbo \noindent{\it Proof of Theorem \ref{invF}. }

For convenience, we write $\cHT =(\ln 2) \cH(T)$. 
$\cHT$ is positive and decreases to $0$. 
So we can write $\cHT = e^{-F_T}$, with $F_T\toi$. Then, 
\eq{
e^{F_T}\phiT \ln \wo{\phiT} + e^{F_T}(1-\phiT)\ln\wo{1-\phiT} + d e^{F_T}\phiT = 1.
}

Note that all three terms on the LHS are positive. 
So as $T\toi$ and $F_T\toi$ we must have $\phiT\to0$.
So we can introduce
\eq{
0\leq \rhoT = (1-\phiT) \phiT e^{F_T} \leq e^{F_T}(1-\phiT)\ln\wo{1-\phiT} \leq 1,
}
since, for $0\leq\phi\leq1$
\eqn{
\label{logineq}
\phi(1-\phi)\leq-(1-\phi)\ln(1-\phi)\leq\phi.
}

Then, multiply through the identity above to find: 
\eq{
\arr{rrcl}{
&&&1-\phiT \\
&&=& (1-\phiT)e^{F_T}\phiT\ln\wo{\phiT}\\
	&&&+e^{F_T}(1-\phiT)^2\ln\wo{1-\phiT}+d\phiT(1-\phiT)e^{F_T}\\
&&= &\rhoT\ln\frac{(1-\phiT)e^{F_T}}{\phiT(1-\phiT)e^{F_T}} 
+ e^{F_T} (1-\phiT)^2\ln\wo{1-\phiT} + d\rhoT\\
&&= &\rhoT\ln\wo{\rhoT} + (1-\phiT)\phiT e^{F_T}\bra{F_T-\ln\wo{1-\phiT}} \\
&&&+ e^{F_T} (1-\phiT)^2\ln\wo{1-\phiT}+d\rhoT\\
&&= &\rhoT\ln\wo{\rhoT} + F_T\rhoT\\
&&& + e^{F_T}(1-\phiT)(1-\phiT)\ln\wo{1-\phiT} + d\rhoT.
}}

Now all terms are $>0$ and $0<\rhoT\leq 1$. So letting $T\toi$ we get  $F_T\toi$ 
and $1-\phiT\ra 1$
and so must have
 $\rhoT\to0$.

Then, as $\Ttoi$, the first term $\rhoT\ln\wo{\rhoT}\to0$, the last term $d\rhoT\to0$ 
and due to (\ref{logineq}), the third term 
\eq{
e^{F_T}(1-\phiT)(1-\phiT)\ln\wo{1-\phiT}&\leq e^{F_T}(1-\phiT)\phiT=\rhoT\to0.
}

So, we are left with, as $\Ttoi$, 
\eq{
	& F_T\rhoT = F_Te^{F_T}\phiT(1-\phiT)\to 1\\
\Ra & F_Te^{F_T} \phiT \to 1\\
\Ra & \phiT \sim \frac{\emu{F_T}}{F_T}=\frac{\cHT}{-\log\cHT}.
}\hfill$\square$



\section*{B. Proofs for Section \ref{aff}}
\label{proaff}

\subsection*{B.1. Proofs for Section \ref{aff} A.}

\noindent{\it Proof of Lemma \ref{RPlb}. }
We have
\eq{
\Dkj(T) 
&= \E_k\left[\ln\frac{L_k(\NoT)}{L_j(\NoT)}\right]\\
	&= \E_k\sqbra{\ssum{r=1}^{N_T} \ln\frac{\lamk(T_r)}{\lamj(T_r)} - \intoT \lamk(t)-\lamj(t)dt}\\
	&=   \E_k\sqbra{\intoT\ln\frac{\lamk(t)}{\lamj(t)}dN_t-\intoT\lamk(t)-\lamj(t)dt}\\
	&= \E_k\sqbra{\intoT\lamk(t)\ln\frac{\lamk(t)}{\lamj(t)} -\lamk(t)+\lamj(t)dt}\\
	&=\intoT\E_k[\lamk(t)\lamkj(t)]dt,
}
where $\lamkj(t) = \frac{\lamj(t)}{\lamk(t)} - 1 - \ln\frac{\lamj(t)}{\lamk(t)}$. 

We replaced the summation with the Lebesgue-Stieltjes integral of the counting measure 
in the third line. And we used the martingale property of $dN_t-\lam(t)dt$ in the second last line
(\cite{Daley03} and \cite{Bremaud20}).

Using the classic hazard relations in Section \ref{prel}, we can rewrite the Janossy density
\eq{
L_k(\NoT) &= \sprod{r=1}^n \lamk(\tr) e^{-\Lamk(T)}\\
&= \sprod{r=1}^n \pk(\tr-\trmw) \Sk(T-t_n),}
with $t_0=0$. 
Further 
\eq{
\lamk(t)\lamkj(t) &= \hk(t-t_n)\bra{\frac{\hj(t-t_n)}{\hk(t-t_n)} - 1 -\ln\frac{\hj(t-t_n)}{\hk(t-t_n)}}\\
	&\equiv\hk(t-t_n)\hkj(t-t_n),
}
We note also that
\eq{
\Bkj(x)&=\pk(x)\hkj(x) = \qkj(x) + \frac{d\Fkj(x)}{dx}.
}

We now find $\E_k[\lamk(t)\lamkj(t)]=\ints \calj_n d\twn$ where
\eq{
\calj_n&= \sprod{r=1}^n \pk(\tr - \trmw)\Sk(t-\tn) \hk(t-\tn)\hkj(t-\tn)\\
&= \sprod{r=1}^n \pk(\tr - \trmw)\pk(t-\tn)\hkj(t-\tn)\\
&= \sprod{r=1}^n \pk(\tr - \trmw)\Bkj(t-\tn).
}
We now integrate backwards from $t_n$, finding
\eq{
	&\int_{t_{n-1}}^t \sprod{r=1}^n \pk(\tr - \trmw)\Bkj(t-\tn) dt_n\\
=& \sprod{r=1}^{n-1}\pk(\tr - \trmw)\int_0^{t-t_{n-1}} \pk(x_n) \Bkj(t-t_{n-1}-x_n)dx_n\\
=& \sprod{r=1}^{n-1}\pk(\tr - \trmw)(\Bkj\star \pk)(t-t_{n-1}),
}
where $(f*g)(t) = \int_0^t f(u)g(t-u)du$ is a convolution. 
Continuing the backwards integration gives
\eq{
\arr{rcll}{
 &\E_k[\lamk(t)\lamkj(t)]  &=& \Bkj(t) + \ssum{n=1}^\infty (\Bkj\star \pk^{(n)})(t)\\
\Ra & \Dbkj(s) &=& \wos \bra{\Bbkj(s) +  \ssum{n=1}^\infty \Bbkj(s) \pbk^n(s)} \\
&&=& \frac{\Bbkj(s)}{s(1-\pbk(s))}.
}}\hfill$\square$

\subsection*{B.2. Proofs for Section \ref{aff} B.} 
Now for Theorem \ref{Ra}, apparently we can get this result from one in 
\cite{Feller88}. 
But \cite{Feller88} does not give existence conditions for 
$\gamkj$ and the integrals in (a) and (b). 
We therefore take a different approach, 
based on the Blackwell's final value theorem (FVT). 

\lem{
\label{Bw}
{Blackwell's Final Value Theorem} \cite{BGT87}. 
Let $g(x)$ have LT $\bar g(s)=\frac{\bar r(s)}{1-\pbar(s)}$, 
where $\pbar(s)$ is the LT of a probability density $p(x)$. 
Then, 
\eq{
\limso s\bar g(s) = c\Ra\limxi g(x) = c= \wo\mu\bar r(0),
}
where $\mu=-\bar p'(0) = \intoi up(u)du$.
}

\noindent
{\it Remark. }
The FVT in traditional engineering applications
assumes the LT is rational. 
Lemma \ref{Bw} is a deeper result, for non-rational LTs, needed for our theory below.\\

\noindent{\it Proof of Theorem \ref{Ra}. }
We have to show
$e^{\gamkj T}\Rkj(T)\to\alpkj$ as $T\toi$.
We use Blackwell's FVT to do this.
Since 
$
\Rbkj(s) = \frac{\Gbkj(s)}{1-\pbkj(s)},
$
then $e^{\gamkj T}\Rkj(T)$ has LT 
\eq{
\Rbkj(s-\gamkj) = \frac{\Gbkj(s-\gamkj)}{1-\pbkj(s-\gamkj)}.
}

To apply the FVT we need to check that $\pbkj(s-\gamkj)$ is
the LT of a probability density i.e. that $\pbkj(-\gamkj)=1$.
However this holds since
$\pbkj(-\gamkj)=\intoi e^{\gamkj x}\pkj(x)dx=\intoi\pkj^*(x)dx=1$.

The FVT now gives
\eq{
\alpkj &= \limso s\Rbkj(s-\gamkj) \\
&= \limso\frac{\Gbkj(s-\gamkj)}{\wos(1-\pbkj(s-\gamkj))}\\
&= \frac{\Gbkj(-\gamkj)}{-\pbkj'(-\gamkj)}
}
By the elementary inequality $2ab\leq a^2+b^2$ we find
\eq{
\pbkj(0) < \intoi \frac12(\pk(x)+\pj(x))dx = 1.
}

By Definition \ref{mgf}, we can find a $\gamkj>0$ such that $\pbkj(-\gamkj)<\infty$. 
Then, since $\pbkj(-\gam)$ is a continuous and increasing function of $\gam$ and 
$\pbkj(-\infty)=\infty$, we can find a unique $0<\gamkj<\infty$, such that 
$\pbkj(-\gamkj) = 1$. \hfill$\square$\\

\bult\mbo\noindent{\it Proof of Lemma \ref{Gp}. }
With A1, we can find $x_k^*,x_j^*$, such that 
$\inf_{x>x_k^*} \hk(x) = z_k>0$ and  $\inf_{x>x_j^*} \hj(x) =z_j>0$. 
Suppose $\Gbkj(-\gamkj)=\infty$. 
Let $x^* = \max\{x_k^*, x_j^*\}$. Then, use the classic hazard relations to find
\eq{
\pbkj(-\gamkj) 
&= g + \int_{x^*}^\infty \sqrt{\hk(x)\hj(x)} e^{\gamkj x} \Gkj(x)dx\\
	&\geq g + \sqrt{z_kz_j} \int_{x^*}^\infty e^{\gamkj x} \Gkj(x)dx =  \infty,
}
where $g = \int_0^{x^*} \pkj^*(x)dx<\infty$. 
This contradicts $\pbkj(-\gamkj)=1$ and the result follows.\hfill$\square$\\

\bult\mbo We are left with the discussion of condition (b): $\intoi x\pkj^*(x)dx<\infty$. 
First, suppose at least one of $\pk$ and $\pj$, say $\pk$, has 
$\intoi e^{\theta x}\pk(x)dx<\infty$
for all real $\theta$, e.g. $\pk(x)\sim e^{-x^2}$. 
Then, it is straightforward to show that 
condition (b) $\intoi x\pkj^*(x)dx=\intoi x\sqrt{\pk(x)\pj(x)}dx<\infty$ is always satisfied, 
since $\pkj^*$ is also regular and thus the first moment exists.

Now suppose there exist $\thek^*, \thej^*>0$, such that the MGFs 
$m_k(\theta) = \intoi e^{\theta x} \pk(x)dx$ and $m_j(\theta) = \intoi e^{\theta x} \pj(x)dx$ 
exist for all $\theta<\thek^*$ and $\theta<\thej^*$, respectively. 
But $m_k(\theta)=\infty, m_j(\theta)=\infty$ for $\theta>\thek^*$ and $\theta>\thej^*$. 
Then we must have $\gamkj\leq\frac{\thek^*+\thej^*}2$. 

For the case where $\gamkj<\frac{\thek^*+\thej^*}2$, clearly $\pkj^*$ is regular 
and thus condition (b) is satisfied. 
For the extreme case $\gamkj=\frac{\thek^*+\thej^*}2$, 
it is possible that the first moment does not exist. 
However, such a case has measure $0$.

\subsection*{B.3. Proofs for Section \ref{aff} C. }
\noindent{\it Proof of Theorem \ref{Ka}. }
Introduce $g(T) = \Dkj(T) - \rhokj T$. 
We need to show that 
as $\stoo$, $s\bar g(s) \to \ln\ckj$. So that, by Blackwell's FVT, 
as $\Ttoi$, 
\eq{
\arr{rrl}{
	&\Dkj(T) - \rhokj T &\to \ln\ckj\\
\equiv & e^{\rhokj T} \Kkj(T)=e^{\rhokj T - \Dkj(T)}&\to\ckj.
}}

So noting that $\Sbk(s)=\frac{1-\pbk(s)}s$, we get
\eq{
s\bar g(s) &= s\Dbkj(s) - \frac{\rhokj}{s}\\
&=\frac{\Bbkj(s)}{1-\pbk(s)} - \frac{\rhokj}{s}\\
	&= \frac{\wos(\qbkj(s)-\rhokj\Sbk(s)) + \Fbkj(s)}{\Sbk(s)}.
}

We then have
\eq{
\limso s\bar g(s) = \frac{\sqbra{\limso\wos\bar\chi(s)} + \Fbkj(0)}{\muk},
}
where $\bar\chi(s) = \qbkj(s) - \rhokj \Sbk(s)$ and (c) 
$\Fbkj(0) = \intoi \Fkj(x)dx<\infty$ is to be proved.

We have $\bar\chi(0) = \qbkj(0) - \frac{\qbkj(0)}{\muk}\muk = 0.$ 
But the LT relations give
\eq{
\limso\bar\chi'(s) &= -\intoi x\qkj(x)dx + \rhokj \intoi x\Sk(x)dx\\
	&= -\intoi x\qkj(x)dx + \frac12\rhokj \intoi x^2\pk(x)dx\\
	&= \qbkj'(0) + \frac12\rhokj (\sigk^2+\muk^2),
}
where $-\qbkj'(0) = \intoi x \pk(x)\ln\frac{\pk(x)}{\pj(x)}dx<\infty$ by A2*,
and $\sigk^2 = \var_k[X]<\infty$ is ensured by existence of the MGF.

So, by L'Hopital's rule, 
\eq{
\arr{rrcl}{
&\limso \wos \chi(s) &=& \qbkj'(0) + \frac12 \rhokj(\sigk^2 + \muk^2)\\
\Ra &\limso s\bar g(s) &=&\ln\ckj,}}
as required.\hfill$\square$\\

We now prove that condition (c) holds under A2*.

\lem{
\label{litl}
Log-integral Tail Lemma. 
Let $\pk(x),\pj(x)$ be densities with survivors $\Sk(x), \Sj(x)$. Then
\eq{
\int_x^\infty \pk(u)\ln\frac{\pk(u)}{\pj(u)}du \geq \Sk(x)\ln\frac{\Sk(x)}{\Sj(x)}.
}

Further,
\eq{
\intoi u\pk(u)\ln\frac{\pk(u)}{\pj(u)}du \geq \intoi \Sk(x)\ln\frac{\Sk(x)}{\Sj(x)}dx.
}
}

\pro{
Let $f,g,F,G$ be positive numbers. In the inequality
$
x\ln x - x+1\geq0, x\geq0,
$
set $x=\frac{f/F}{g/G}$. Then, 
\eq{
\arr{rrcl}{
	&0&\leq&\frac{f/F}{g/G}\ln \frac{f/F}{g/G} - \frac{f/F}{g/G}+1\\
\Ra &0&\leq&\frac fF \ln \frac{f/F}{g/G} - \frac fF +\frac gG\\
\Ra &0&\leq&\frac fF \ln \frac fF - \frac fF \ln\frac gG - \frac fF +\frac gG.
}}

Now set $f=\pk(u), g=\pj(u), F=\Sk(u), G=\Sj(u)$ and integrate $u$ from $x$ to $\infty$ to find
\eq{
\wo{\Sk(x)}\int_x^\infty \pk(u)\ln\frac{\pk(u)}{\pj(u)}du - \ln\frac{\Sk(u)}{\Sj(u)}\geq0, 
}
which delivers the first result. 

Then, integrating through the first inequality gives
\eq{
\intoi \int_x^\infty \pk(u)\ln\frac{\pk(u)}{\pj(u)}dudx \geq \intoi \Sk(x)\ln\frac{\Sk(x)}{\Sj(x)}dx.
}

Changing the order of integration on the left-hand side gives the quoted second result.
Thus  A2* $\Ra$ (c).}

\section*{C. Proofs for Section \ref{main}}
\label{promain}

\noindent{\it Proof of Theorem \ref{ab}. }
The first two parts are just the results in Theorem \ref{Ra} and \ref{Ka}.
By the equivalence $\log(1+f)\sim \wo{\ln2}f$ as $f\to0$, we have
\eq{
(\ln2)\cH_{R}(T) &\sim \skwM\ssum{j\neq k} \pi_j \Rkj(T)\sim \alpha^* e^{-\gam^* T}\\
(\ln2)\cH_K(T) &\sim\skwM\ssum{j\neq k} \pi_j \Kkj(T)\sim c^* e^{-\rho^* T}, 
}
where $\alpha^* e^{-\gam^* T}$ and $c^* e^{-\rho^* T}$ are the 
most slowly decaying terms among $\Rkj(T)$ and $\Kkj(T)$, respectively.

The quoted bounds are from Lemma \ref{ulb2}(c).  For the upper bound, 
\eq{
\frac12 \cH_R(T) \sim \wo{2\ln2}\alpha^* e^{-\gam^* T}.
}

For the lower bound, we have $\Kkj(T)\sim\ckj e^{-\rhokj T}\to0$ as $\Ttoi$. 
Then obviously $\cH_K(T)\to0$. 
We then apply Theorem \ref{invF} to find
\eq{
\phiT\sim\frac{\calhk(T)}{-\log\calhk(T)}\sim\frac{\frac{c^*}{\ln2} e^{-\rho^* T}}{\frac{\rho^*}{\ln2} T} = \frac{c^*}{\rho^*}\wo{T} e^{-\rho^* T}.
}\hfill$\square$

\bibliographystyle{plain}
\bibliography{vs-B,xr}  
\end{document}